\documentclass{amsart}
\usepackage{amssymb}
\usepackage{amscd,enumerate,amsfonts,calc,amsmath,verbatim,xypic}

\newcommand{\ann}{\operatorname{ann}}

\newcommand{\Supp}{\operatorname{Supp}}
\newcommand{\Ass}{\operatorname{Ass}}
\newcommand{\Spec}{\operatorname{Spec}}

\newcommand{\gr}{\operatorname{gr}}

\newcommand{\pd}{\mathrm{projdim}}
\newcommand{\depth}{\operatorname{depth}}
\newcommand{\codepth}{\operatorname{codepth}}

\newcommand{\height}{\operatorname{ht}}

\newcommand{\ov}{\overline}
\newcommand{\ovM}{\,\overline{\phantom{O}\!\!\!\!\!\!M}}
 \newcommand{\ovN}{\,\overline{\phantom{O}\!\!\!\!\!\!N}}

\newcommand{\wM}{\,\widetilde{\phantom{O}\!\!\!\!\!\!M}}
\newcommand{\hM}{\,\widehat{\phantom{A}\!\!\!\!\!\!M}}

\newcommand{\col}{\colon}

\newcommand{\m}{{\mathfrak m}}
\newcommand{\n}{{\mathfrak n}}
\newcommand{\p}{{\mathfrak p}}
\newcommand{\q}{{\mathfrak q}}
\newcommand{\fb}{{\mathfrak b}}

\theoremstyle{remark}

\swapnumbers

\theoremstyle{plain}
\newtheorem{theorem}{Theorem}[section]

\newtheorem{proposition}[theorem]{Proposition}
\newtheorem{lemma}[theorem]{Lemma}
\newtheorem{corollary}[theorem]{Corollary}
\newtheorem*{proposition*}{Proposition}
\newtheorem*{lemma*}{Lemma}

\theoremstyle{definition}
\newtheorem{definition}[theorem]{Definition}

\newtheorem*{chunk*}{}

\theoremstyle{remark}
\newtheorem*{Claim1}{Claim 1}
\newtheorem*{Claim2}{Claim 2}

\numberwithin{theorem}{subsection} \numberwithin{equation}{theorem}

\numberwithin{subchunk}{theorem}

\def\zz{\mathbb{Z}}
\def\qq{\mathbb{Q}}
\def\nn{\mathbb{N}}

\begin{document}
\title[Open loci]
{Open loci of graded modules}
\author[C.~Rotthaus]{Christel Rotthaus}
\address{Department of Mathematics, Michigan State University, East Lansing, MI 48824}
\email{rotthaus@math.msu.edu}

\author[L.~M.~\c Sega]{Liana M.~\c Sega}
\address{Department of Mathematics, Michigan State University, East Lansing, MI 48824}
\email{lsega@math.msu.edu}
\begin{abstract}
Let $A=\oplus_{i\in \nn}A_i$ be an excellent homogeneous Noetherian graded ring and let $M=\oplus_{n\in \zz}M_n$
be a finitely generated graded $A$-module. We consider $M$ as a module over $A_0$ and show that the $(S_k)$-loci
of $M$ are open in $\Spec(A_0)$. In particular, the Cohen-Macaulay locus $U^0_{CM}=\{\p\in \Spec(A_0) \mid M_\p
\; \mbox{is Cohen-Macaulay}\}$ is an open subset of $\Spec(A_0)$. We also show that the $(S_k)$-loci on the
homogeneous parts $M_n$ of $M$ are eventually stable. As an application we obtain that for a finitely generated
Cohen-Macaulay module $M$ over an excellent ring $A$ and for an ideal $I\subseteq A$ which is not contained in
any minimal prime of $M$ the $(S_k)$-loci for the modules $M/I^nM$ are eventually stable.

\end{abstract}
\date{\today}
\maketitle
\section*{introduction}
A well-known theorem of Grothendieck states that if $M$ is a finitely generated module over an excellent
Noetherian ring then for all $k\in \nn$ the $(S_k)$-locus of $M$:
$$U_{S_k}(M) =\{\p\in\Spec(A)\mid M_\p \; \mbox{satisfies}\; (S_k)\}$$
is an open subset of $\Spec(A)$. As usual, $(S_k)$ denotes the Serre condition, that is, $M_\p$ satisfies
$(S_k)$ if for all $\q\in\Spec(A)$ with $\q\subseteq \p$ it holds that:
$$ \mbox{depth}_{A_\q}(M_\q)\ge \mbox{min}(k,\mbox{dim}(M_\q)).$$
It also follows that for such modules $M$ the Cohen-Macaulay locus:
$$U_{CM}(M)= \{\p\in \Spec(A) \mid M_\p \; \mbox{is Cohen-Macaulay}\}$$
is an open subset of $\Spec(A)$.

Let $A=\bigoplus_{n\ge 0}A_n$ be Noetherian graded excellent homogeneous ring and $M=\bigoplus_{i\in \zz}M_i$ a
finitely generated graded $A$-module. Considered as a module over the base ring $A_0$, $M$ is a direct sum of
finitely generated $A_0$-modules. Moreover, if the base ring $A_0$ is local the standard notion of depth is
meaningful for the $A_0$-module $M$ and we may consider its $(S_k)$-loci:
$$U^0_{S_k}(M) = \{\p\in \Spec(A_0)\mid M_\p \; \mbox{satisfies}\;
S_k\}$$ where $M_\p$ denotes the localization of $M$ at the multiplicative set $A_0\smallsetminus\p$. In this
paper we prove that under these assumptions the $(S_k)$-loci of the $A_0$-module $M$ are open subsets of
$\Spec(A_0)$. In particular, the Cohen-Macaulay locus of $M$ (as an $A_0$-module):
$$U^0_{CM}(M) = \{\p\in \Spec(A_0)\mid M_\p \; \mbox{is Cohen-Macaulay}\}$$
is an open subset of $\Spec(A_0)$.

The proof follows the main ideas of Grothendieck's proof. It is, however, not merely a copy of the proof in EGA
and requires a number of modifications. For the benefit of the reader we have included complete proofs of the
results. Our proof is based on the following two observations: First, if $A$ is a polynomial ring over the base
ring $A_0$, then every graded resolution of $M$ by finitely generated graded free $A$-modules provides a free
resolution of the $A_0$-module $M$ which is finitely generated on the homogeneous parts. The second is a result
by Hochster and Roberts which states for the $A$-module $M$ that there is an element $a\in A_0-(0)$ so that
$M_a$ is a free $(A_0)_a$-module provided that the ring $A_0$ is a domain.

The paper is organized as follows:

The first section contains basic facts about graded rings and modules which are relevant for the rest of the
paper. As a main result we obtain that the Auslander-Buchsbaum formula holds for the $A_0$-module $M$.

The second section shows that the codepth-loci of $M$ are open in $\Spec(A_0)$. This is the main step in proving
the openness of the $(S_k)$-loci which we present in the next section.

In Section $4$ we consider the homogeneous parts of the graded module $M$. We show that the codepth-loci and
$(S_k)$-loci of the homogeneous parts of $M$ are eventually stable. This is applied in the last section to the
case of a finitely generated module $M$ over an excellent Noetherian ring $A$. If $I\subseteq A$ is an ideal we
recover a well-known result by Kodiyalam \cite{Ko}, namely that for $k\ge k_0$:
$$ \mbox{depth}(M/I^kM) = \mbox{depth}(M/I^{k_0}M).$$
We also show that if $M$ is a Cohen-Macaulay module over $A$ and if $I\subseteq A$ is not contained in a minimal
prime of $M$, then the codepth- and $(S_k)$-loci of $M/I^nM$ are eventually stable.

\section{basic facts}

In this paper we assume that $A = \bigoplus_{i\in \nn} A_i$ is a Noetherian homogeneous graded ring and that $M
=\bigoplus_{i\in \zz} M_i$ is a finitely generated $A$-module. As usual, we let $A_{+}$ denote the irrelevant
ideal of $A$, that is, $A_+=\bigoplus_{i\ge 1}A_i$.

If $\p\in \Spec(A_0)$ is a prime ideal of $A_0$, then $M_\p$ denotes the localization $S^{-1}M$ where $S
=A_0\smallsetminus\p$. Note that $M_\p$ is a graded module over the graded ring $A_\p$.

Our goal is to show that if $A$ is excellent then the codepth-loci and the $(S_k)$-loci of $M$, considered as a
module over the base ring $A_0$, are open subsets of $\mbox{Spec}(A_0)$.

\subsection {General remarks}
We begin our investigation with some well known facts about graded modules. Since these results are frequently
used throughout the paper we include them together with their (short) proofs in this introductory section.


\begin{lemma}
\label{ann-equal} There exists an integer $t$ so that $\ann_{A_0}(M_t)=\ann_{A_0}(M_k)$ for all $k\ge t$.

\end{lemma}

\begin{proof}
  For all $k\in \zz$ set $J_k = \ann_{A_0}(M_k)$ . Since $A$ is
  homogeneous and $M$ is a finitely generated $A$-module, there exists
  $s_0\in \zz$ such that
$$A_1M_k = M_{k+1}\quad\text{for all}\quad k\ge t_0\,.$$
We conclude $ J_k\subseteq J_{k+1}$ for all $k\ge t_0$.  Since $A_0$ is Noetherian, there exists then $s\ge t_0$
so that $J_k = J_t$ for all $k\ge s$.
\end{proof}

\begin{lemma}
\label{two-functions} The following two functions are well-defined and surjective:
\begin{enumerate}[\quad\rm(1)]
\item The function $\varphi\col\Supp_A(M)\to\Supp_{A_0}(M)$ defined by
  \,$\varphi(P) = P\cap A_0$.
\item The function $\psi\col \Ass_A(M)\to\Ass_{A_0}(M)$ defined by
  \,$\psi(P) =P\cap A_0$.
\end{enumerate}
\end{lemma}

\begin{proof}
  $(1)$ If $P\in \mbox{Supp}_A(M)$, then $M_P\ne 0$ and in particular
  $M_\p\ne 0$, where $\p = P\cap A_0$. This shows that $\varphi$ is
  well defined. Let $\p\in \mbox{Supp}_{A_0}(M)$, then
\[
M_\p = \bigoplus_{i\in \zz} (M_i)_\p\ne 0
\]
and we may consider $M_\p$ as a graded module over the graded ring $A_\p$. Note that $A_\p$ is a $^*$local ring
with unique graded maximal ideal $\m=\p(A_0)_\p\bigoplus(A_{+})_\p$. Since all minimal primes of
$\Supp_{A_\p}(M_\p)$ are graded, $\m\in \Supp_{A_\p}(M_\p)$. Thus there is a prime $P\in \Supp_A(M)$ with $P\cap
A_0=\p$.

$($2$)$ If $P\in \Ass_A(M)$, then there exists $y\in M$ so that $\ann_A(y)=P$. Thus $\ann_{A_0}(y)=P\cap A_0
=\p$ and $\p\in \Ass_{A_0}(M)$. Conversely, let $\p\in \Ass_{A_0}(M)$. Consider again the graded $A_\p$-module
$M_\p$.  There exists $z\in M_\p$ so that $\ann_{(A_0)_\p}(z) = \p(A_0)_\p$, and therefore
$$\p(A_0)_\p\subseteq\bigcup_{Q\in \Ass_{A_\p}(M_\p)} Q$$
Since $M_\p$ is a finitely generated $A_\p$-module, there exists $Q\in \Ass_{A_\p}(M_\p)$ with
$\p(A_0)_\p\subseteq Q$. Since $A_\p$ is $^*$local with unique graded maximal ideal $\p(A_0)_\p \oplus
(A_{+})_\p$ we obtain $Q\cap(A_0)_\p = \p(A_0)_\p$ and a preimage $P\in\Spec(A)$ of $Q$ is an associated prime
of the $A$-module $M$, with $P\cap A_0 = \p$.
\end{proof}
Lemma \ref{two-functions} shows in particular that $M$ as an $A_0$-module has a finite set of associated primes.

\begin{lemma}
\label{ann} Let $A$ and $M$ be as above and set $I=\ann_{A_0}(M)$. For any $\p\in\Spec(A_0)$ the following hold:
\begin{enumerate}[\quad\rm(1)]
\item If $M_\p=0$, then there is an element $a\in A_0\smallsetminus\p$
  with $M_a=0$.
\item $\ann_{(A_0)_\p}(M_\p) = I(A_0)_\p$
\end{enumerate}
\end{lemma}

\begin{proof}
  $($1$)$ This is a basic fact about Noetherian modules using that $M$ is
  a finitely generated module over $A$ and $A_0\smallsetminus\p$ is a
  multiplicative subset of $A$.

  $($2$)$ Obviously, $I(A_0)_\p\subseteq \ann_{(A_0)_\p}(M_\p)$.  Let
  $x\in \ann_{(A_0)_\p}(M_\p)$ with $x = \frac{b}{s}$ where $b\in A_0$
  and $s\in A_0\smallsetminus\p$. Assume that $m_1,\hdots,m_r$ is a
  system of generators of the $A$-module $M$. Since $x\frac{m_i}{1}=0$
  for all $1\le i\le r$ there is an element $t\in
  A_0\smallsetminus\p$ with $tbm_i = 0$ for all $1\le i\le r$. We have
  that $tb\in I$ and hence $x=\frac{b}{s}\in I(A_0)_\p$.
\end{proof}

\subsection {The Auslander-Buchsbaum formula}

Let $A=\bigoplus_{i\ge 0}A_i$ be a graded Noetherian homogeneous ring with $(A_0,\m_0)$ local and let
$M=\bigoplus_{i\in \zz} M_i$ be a finitely generated $A$-module. Since $M$ is (in general) not finitely
generated as an $A_0$-module, we need to verify that the classical definition of $A_0$-depth works in the case
of a finitely generated graded module. First note that an element $z\in \m_0$ is regular on $M$ if and only if
$z$ is regular on $M_i$ for all $i\in \zz$ with $M_i\ne 0$. Let $x_1,\hdots,x_s\in \m_0$ and $y_1,\hdots,y_t\in
\m_0$ be two maximal regular $M$-sequences (as an $A_0$-module). Then for all $i\in \zz$ with $M_i\ne 0$ the two
sequences are regular on the $A_0$-module $M_i$ and the sets
\begin{gather*}\Ass_{A_0}(M/(x_1,\hdots,x_s)M) = \bigcup_{i\in \zz}
\Ass_{A_0}(M_i/(x_1,\hdots,x_s)M_i)\\
\Ass_{A_0}(M/(y_1,\hdots,y_t)M) = \bigcup_{i\in \zz} \Ass_{A_0}(M_i/(y_1,\hdots,y_t)M_i)
\end{gather*}
are finite by Lemma \ref{two-functions}. The maximality of the first sequences yields
 that there is an
$i\in \zz$ with $M_i\ne 0$ and $\m_0\in \Ass_{A_0}(M_i/(x_1,\hdots,x_s)M_i)$. Since the second sequence is also
regular on $M_i$ we have that $t\le s$. A similar argument shows that $s\le t$ and we obtain that two maximal
regular sequences on $M$ have the same length. Therefore the classical definition of depth is efficient and we
put:

\begin{definition}
\label{depth-defn}
  Let $A$ and $M$ be as above with $(A_0,\m_0)$ local. We define the
  {\it depth} of $M$ as $A_0$-module to be the number:
  $$\depth_{A_0}(M) := \sup\{ n\in \nn \mid \exists \,\,\,
  \mbox{an $M$-sequence of lenght $n$} \}\,.$$
\end{definition}

\begin{lemma}
\label{depth-projdim} Let $A$ and $M$ be as above and assume that $(A_0, \m_0)$ is local. Then:
\begin{enumerate}[\quad\rm(1)] \item $\dim_{A_0}(M) =\sup\{\dim_{A_0}(M_i) \mid i\in \zz \}$
\item $\depth_{A_0}(M) = \inf\{\depth_{A_0}(M_i)\mid i\in \zz \;
  \mathrm{with} \; M_i \ne 0 \}$
\item $\pd_{A_0}(M) = \sup\{\pd_{A_0}(M_i)\mid i\in \zz \}$
\end{enumerate}
\end{lemma}

\begin{proof}
$($1$)$ By Lemma \ref{ann-equal} there is an integer $s\in \zz$ so that $\mbox{ann}_{A_0}(M_k) =
\mbox{ann}_{A_0}(M_s)$ for all $k\ge s$. In particular, for all $k\ge s$: $\mbox{dim}_{A_0}(M_k) =
\mbox{dim}_{A_0}(M_s)$ and
$$\mbox{dim}_{A_0}(M) = \mbox{dim}_{A_0}(M_r \oplus M_{r-1} \oplus \hdots \oplus M_{s-1}\oplus M_s)$$ where $r\in \zz$ is the
smallest integer $j$ with $M_j\ne 0$. The dimension of a finite direct sum of $A_0$-modules is the maximum of
the dimensions of its summands.

  $($2$)$ If $r_1,\hdots,r_s\in A_0$ is a regular sequence on $M$,
  then $r_1,\hdots,r_s$ is a regular sequence on $M_i$ for all $i\in
  \zz$ with $M_i\ne 0$. Thus $ \depth_{A_0}(M) \le \depth_{A_0}(M_i)$
  for all $i\in \zz$ with $M_i\ne 0$ and hence
  $$\depth_{A_0}(M) \le \inf\{\depth_{A_0}(M_i)\mid i\in \zz \;
  \mbox{with} \; M_i \ne 0 \}\,.$$
  In order to show the other
  inequality we proceed by induction on $t =\depth_{A_0}(M)$.

Note that by Lemma \ref{ann}, $\Ass_{A_0}(M)$ is a finite set.

If $t=0$ then $\m_0\in \Ass_{A_0}(M)$ and there is an $i\in \zz$ so that $\m_0\in \Ass_{A_0}(M_i)$. Thus:
$$\inf\{\depth_{A_0}(M_i)\mid i\in \zz \; \mbox{with} \; M_i \ne
0 \} = 0\,.$$ Now assume that $t= \depth_{A_0}(M)> 0$. This implies that
$$\bigcup_{\p\in \Ass_{A_0}(M)} \p \ne \m_0\,.$$
Consider an element
$$r\in \m_0\smallsetminus\bigcup_{\p\in \Ass_{A_0}(M)} \p$$
Since $r $ is regular on $M$, and therefore is regular on $M_i$ for all $i\in \zz$ with $M_i\ne 0$, we obtain
$$\depth_{A_0}(M/rM) = \depth_{A_0}(M)-1$$
and for all $i\in \zz$ with $M_i\ne 0$:
$$\depth_{A_0}(M_i/rM_i) =\depth_{A_0}(M_i) -1\,.$$

By induction hypothesis:
$$\depth_{A_0}(M/rM)= \inf\{ \depth_{A_0}(M_i/rM_i)\mid i\in \zz \;
\mbox{and} \; M_i/rM_i \ne 0\}\,.$$ The assertion follows.

$($3$)$ For all $i\in \zz$ let $F^{(i)}_{\bullet}$ be a finite free resolution of $M_i$. Then
$$ F_{\bullet} = \bigoplus_{i\in \zz} F^{(i)}_{\bullet} $$
is a free resolution of the $A_0$-module $M$ yielding:
$$
\pd_{A_0}(M) \le \sup\{ \pd_{A_0}(M_i)\mid i\in \zz \}\,.$$ In order to show the other inequality, assume that
$\pd_{A_0}(M) =r$ and consider for all $i\in \zz$ the $r$th-syzygy $T^{(i)}_r$ of $M_i$ and the exact sequence:
$$
0\longrightarrow T_r^{(i)} \longrightarrow F_{r-1}^{(i)} \longrightarrow \hdots \longrightarrow F_0^{(i)}
\longrightarrow M_i \longrightarrow 0\,.$$ By taking direct sums we see that
$$
\bigoplus_{i\in \zz} T_r^{(i)}$$ is an $r$th-syzygy of $M$ and thus projective. Therefore every $T_r^{(i)}$ is a
projective finitely generated $A_0$-module. Since $A_0$ is a local Noetherian ring every, $T_r^{(i)}$ is a free
$A_0$-module and thus for all $i\in \zz$:
$$ \pd_{A_0}(M_i)\le r\,.$$
This shows $($3$)$.
\end{proof}

\begin{proposition}
\label{AB-formula}
  Let $A$ and $M$ be as above with $(A_0,\m_0)$ a local ring. Then the
  Auslander-Buchsbaum formula holds for $M$ as an $A_0$-module. That
  is, if $\pd_{A_0}(M)$ is finite, then:
  $$\depth_{A_0}(M) + \pd_{A_0}(M) = \depth(A_0)\,.$$
\end{proposition}

\begin{proof}
  Let $\pd_{A_0}(M)=r< \infty$, then by Lemma \ref{depth-projdim}(2)
  there is an $i\in\zz$ with $\pd_{A_0}(M) =\pd_{A_0}(M_i)$ and for
  all $j\in \zz$:
  $$
  \pd_{A_0}(M_j)\le r.$$
  The Auslander-Buchsbaum formula
  holds for finitely generated $A_0$-modules:
$$\depth_{A_0}(M_j)+\pd_{A_0}(M_j) = \depth_{A_0}(A_0)\quad\text{for
all}\quad j\in\zz$$ and therefore:
$$\depth_{A_0}(M_j) \ge \depth_{A_0}(M_i)\quad\text{for all}\quad
j\in\zz$$ Using Lemma \ref{depth-projdim}(1), we conclude $\depth_{A_0}(M)=\depth_{A_0}(M_i)$. The
Auslander-Buchsbaum formula for $M_i$ gives then the desired formula.
\end{proof}

\subsection{An example}

If $M$ is just any infinite sum of finitely generated modules over a Noetherian local ring, then Lemma
\ref{depth-projdim}(3) still holds (with the same proof), but Lemma \ref{depth-projdim}(2) and Proposition
\ref{AB-formula} may fail. Here is an example:

Let $A_0 = \qq[x,y]_{(x,y)}$ and let $\{\p_i\}_{i\in \nn}$ be an enumeration of the height one prime ideals of
$A_0$. Put:
$$ M = \bigoplus_{i\in \nn} A_0/\p_i$$
and set $M_i = A_0/\p_i$. Obviously for all $i\in \nn$:
$$
\depth_{A_0}(M_i) = \pd_{A_0}(M_i)= 1\,.$$ In particular, $\pd_{A_0}(M) = 1$. But $\depth_{A_0}(M) = 0$, since
every non-unit $r\in A_0$ is contained in some prime $\p_i$, hence it is a zerodivisor on $M_i$ and on $M$. This
shows that the Auslander-Buchsbaum formula does not hold for $M$. (The key is that in the proof of Lemma
\ref{depth-projdim}(2) we make use of the fact that the set of associated primes of $M$ is finite, while in this
example $\Ass_{A_0}(M) = \mbox{Spec}^1(A_0)$ is infinite and $\cup_{i\in\mathbb N}\p_i=(x,y)$, the maximal ideal
of $A_0$.)

\section{Openness of the codepth locus}

Throughout this section we assume that $A=\bigoplus_{i\in \nn_0}A_i$ is a graded Noetherian homogeneous ring and
that $M =\bigoplus_{i\in \zz}M_i$ is a finitely generated $A$-module. Our aim is to generalize and/or modify
existing theorems for finitely generated modules over Noetherian rings to the graded case where the module $M$
is considered a module over the base ring $A_0$. We begin with a result on the flat locus of the $A_0$-module
$M$.
\subsection{The flat locus of $M$}\label{flat-locus} Our first result is a modification of \cite[Theorem 24.3]{Ma}.
The proof follows the proof in Matsumura's book. A key observation is that for a finitely generated graded
module $M$ the localizations $M_{\p}$ are $I$-adically separated for every ideal $I\subseteq (A_0)_{\p}$.
\begin{proposition*} Let $A$ and $M$ be as above. The flat locus of $M$ as an $A_0$-module:
$$ U^0(M) = \{ \p\in \Spec(A_0) \; | \; M_{\p} \; \mbox{is flat over} \; A_0\}$$
is open in $\Spec(A_0).$
\end{proposition*}
\begin{proof} According to Nagata's criterion on the openness of loci \cite[Theorem 23.2]{Ma} we have to
show:
\begin{enumerate}[(a)] \item If $\p,\q\in \Spec(A_0)$ with $\p\in U^0(M)$ and $\q\subseteq \p$ then $\q\in
U^0(M)$. \item If $\p\in U^0(M)$ then $U^0(M)$ contains a nonempty open subset of $V^0(\p) = \{\n\in
\Spec(A_0)\; |\; \p\subseteq \n\}$. \end{enumerate} $($a$)$ is trivial. Let $\p\in U^0(M)$, that is, assume that
$M_{\p}$ is flat over $A_0$. Set $\bar A_0= A_0/\p$. By \cite[Theorem 22.3]{Ma} for every $\q\in V^0(\p)$ the
module $M_{\q}$ is flat over $A_0$ if and only if $(M/\p M)_{\q}$ is flat over $\bar A_0$ and
$\mbox{Tor}^{A_0}_1(M_{\q},\bar A_0) =0$. A similar argument as in the proof of \cite[Theorem 23.2]{Ma} shows
that $\mbox{Tor}^{A_0}_1(M,\bar A_0)$ is a finitely generated module over $A$. Therefore there is an element
$a\in A_0\smallsetminus \p$ so that $(\mbox{Tor}^{A_0}_1(M,\bar A_0))_a = 0$. By applying \cite[Theorem
24.1]{Ma} to the $\bar A_0$-module $M/\p M$ we obtain an element $b\in A_0\smallsetminus \p$ so that $(M/\p
M)_b$ is a free $(\bar A_0)_b$-module. Set $D^0_{ab} = \{\q\in \Spec(A_0)\; |\; ab\notin \q\}$, then for all
$\q\in V^0(\p)\cap D^0_{ab}$ we have that $\mbox{Tor}^{A_0}_1(M_{\q},\bar A_0)=0$ and that $(M/\p M)_{\q}$ is
flat over $(\bar A_0)_{\q}$. Thus by \cite[Theorem 22.3]{Ma} the module $M_{\q}$ is flat over $(A_0)_{\q}$ and
$M_{\q}$ is flat over $A_0$.
\end{proof}

\subsection{A proposition by Auslander}
\label{Auslander-prop} As before let $A$ be a Noetherian graded homogeneous ring and let $M$ be a finitely
generated $A$-module. The following Proposition is an extension of a proposition in EGA \cite[(6.11.1) and
(6.11.2)]{Gr} to the (not finitely generated) $A_0$-module $M$.

\begin{proposition*} The function $ \gamma:\Spec(A_0)\longrightarrow \nn$
defined by
$$\gamma(\p) = \pd_{(A_0)_\p}(M_\p)\quad\text{for
all}\quad\p\in\Spec(A_0)$$ is upper semicontinuous. That is, for all $n\in \nn$ the set
$$U^0_n(M) = \{\p\in\Spec(A_0)\mid \pd_{(A_0)_\p}(M_\p)\le n\}$$ is open
in $\mbox{Spec}(A_0)$.
\end{proposition*}
\begin{proof}
  Note that the ring $A$ is the homomorphic image of the polynomial
  ring $B =A_0[x_1,\hdots,x_t]$, and that, with the standard grading
  on the polynomial ring $B$, the graded $B$-module $M$ is finitely
  generated. We may replace $A$ by $B$ and assume that $A$ is a graded
  polynomial ring over $A_0$. Let $\p\in\Spec(A_0)$ with
  $\pd_{(A_0)_\p}(M_\p)\le n$.

  Consider a graded finitely generated free resolution of the
  $A$-module $M$:
  $$
  F_n \xrightarrow{\varphi_n} F_{n-1}
  \xrightarrow{\varphi_{n-1}}\cdots
  \xrightarrow{\varphi_1}F_1
  \xrightarrow {\varphi_0}M \to 0$$
  where
  the $F_i$ are finitely generated graded free $A$-modules and the
  $\varphi_i$ are homogeneous $A$-linear maps. Let $T$ be the
  $n$th syzygy of $M$, yielding an exact sequence of graded
  $A$-modules:
  $$
  (^*) \quad 0 \to T \xrightarrow{\,\delta\,} F_{n-1}
  \xrightarrow{\varphi_{n-1}}\hdots \xrightarrow{\varphi_1} F_1
  \xrightarrow{\varphi_0} M\to 0\,.$$
  Since all the homogeneous parts
  of $F_i$ are free $A_0$-modules and since $T$ is a graded $A$-module
  we obtain for all $k\in \zz$ an exact sequence of $A_0$-modules:
$$0 \to T_k \xrightarrow{(\delta)_k} (F_{n-1})_k
\xrightarrow{(\varphi_{n-1})_k}\hdots \xrightarrow{(\varphi_1)_k} (F_1)_k \xrightarrow{(\varphi_0)_k} M_k\to 0$$
with $(F_i)_k$ a finitely generated free $A_0$-module. Therefore by considering $(^*)$ as an exact sequence of
$A_0$-modules we obtain that every module $F_i$ is free over $A_0$ and $T$ is an $n$th syzygy of the
$A_0$-module $M$. Localization at $\p$ yields exact sequences:
$$0 \to T_\p \xrightarrow{\delta_\p} (F_{n-1})_\p
\xrightarrow{(\varphi_{n-1})_\p}\hdots \xrightarrow{(\varphi_1)_\p} (F_1)_\p \xrightarrow{(\varphi_0)_\p} M_\p
\to 0\,.$$ Since $\pd_{(A_0)p}(M_\p)\le n$ it follows that $T_\p$ is a projective $(A_0)_\p$-module. Therefore
$T_\p$ is a free $(A_0)_\p$-module. Since $T$ is a finitely generated graded $A$-module it follows from
Proposition \ref{flat-locus} that the set
$$ U^0(T) = \{ \q\in\Spec(A_0)\mid T_\q \quad \mbox{is a flat over} \quad
(A_0)_\q\}$$ is an open subset of $\mbox{Spec}(A_0)$. Since $T$ is a finitely generated graded $A$-module:
$$ T = \bigoplus_{i\in \zz} T_i$$
we have for $\q\in\Spec(A_0)$
$$ T_\q = \bigoplus_{i\in \zz} (T_i)_\q\,.$$
If $T_\q$ is flat over $(A_0)_\q$ then, by \cite[chapter 1, \S 2.3, Proposition 2]{Bo}, for all $i\in \zz$,
$(T_i)_\q$ is flat over $(A_0)_\q$. Since every $(T_i)_\q$ is a finitely generated $(A_0)_\q$-module, each
$(T_i)_\q$ is a free $(A_0)_\q$-module and
$$ U^0(T) = \{  \q\in\Spec(A_0)\mid T_\q \quad \mbox{is a free over} \quad
(A_0)_\q\}\,.$$ This shows that $\p\in U^0(T)$ and
$$ U^0(T) \subseteq \{\q\in\Spec(A_0)\mid \pd_{(A_0)_\q}(M_\q)\le n
\}\,.$$ The set $\{\q\in\Spec(A_0)\mid \pd_{(A_0)_\q}(M_\q)\le n\}$ is thus open in $\mbox{Spec}(A_0)$.
\end{proof}

\subsection{A dimension formula}
\label{dim-formula}
\begin{proposition*} Let $A$ and $M$ be as above. Assume that $A_0$ is
catenary and let $\p$ be a prime ideal in $A_0$ with $\p\in \Supp_{A_0}(M)$. Then there is an open subset $U$ in
$\Spec(A_0)$ such that  $\p\in U$ and for all $\q\in U \cap V^0(\p)$ we have:
$$ \dim(M_\q) = \dim(M_\p) + \dim((A_0/\p)_\q)\,.$$
\end{proposition*}

\noindent{\it Proof.} Set $S = A_0/\ann_{A_0}(M)$ and choose an element $a\in S\smallsetminus\p$ so that the
following equality on the set of minimal primes holds:
  $$
  \mbox{Min}(S_\p) = \mbox{Min}(S_a)\,.$$
  Assume that $\mbox{dim}(M_{\p}) =\height(\p S) = t$
  and choose elements $y_1,y_2,\hdots,y_t\in S$ so that:
$$y_1 \quad \mbox{not in a minimal prime of} \quad S_\p$$
$$y_2 \quad \mbox{not in a minimal prime of} \quad y_1S_\p$$
$$\hdots $$
$$y_t \quad \mbox{not in a minimal prime of} \quad (y_1,\hdots,
y_{t-1})S_\p.$$ Then there is an element $b\in S\smallsetminus\p$ so that:
$$y_1 \quad \mbox{not in a minimal prime of} \quad S_b$$
$$y_2 \quad \mbox{not in a minimal prime of} \quad y_1S_b$$
$$\hdots $$
$$y_t \quad \mbox{not in a minimal prime of} \quad (y_1,\hdots,
y_{t-1})S_b\,.$$ Let $a,b$ also denote preimages of $a$ and $b$ in $A_0$ and put $U = D_{ab} =
\{\q\in\Spec(A_0)\mid ab\notin\q\}$. Then for every $\q\in U\cap V^0(\p)$ the elements $y_1,\hdots,y_t$ extend
to a system of parameters of $S_\q$. Since $S_\p$ and $S_\q$ have the same set of minimal primes and since $S$
is catenary we obtain that:
$$ \dim(S_\q) = \dim(S_\p) + \dim((S/\p)_\q)\,.$$
This is the same as:
\begin{xxalignat}{3}
&\ &\dim(M_\q)&=\dim(M_\p) + \dim((A_0/\p)_\q)\,.&&\square
\end{xxalignat}

\subsection{The special case of $A_0$ regular}
\label{regular-case} Let $(R,\m)$ be a local Noetherian ring and $M$ an $R$-module. Then we define:
$$\codepth_R(M) := \dim_R(M) - \depth_R(M)\,.$$
As usual the depth of the zero module is defined to be $\infty$, and the dimension of the zero module is
$-\infty$, implying that the codepth of the zero module is $-\infty$.

The following proposition extends a result by Auslander \cite[(6.11.2)]{Gr} to the graded case.

\begin{proposition*} Let $A$ and $M$ be as above and assume  that $A_0$ is
a homomorphic image of a regular ring. The function $\varphi\col\Spec(A_0)\longrightarrow \nn$ defined by
\[\varphi(\p) = \codepth_{(A_0)_\p}(M_\p)\quad\text{for all}\quad\p\in\Spec(A_0)
\]
is upper semicontinuous, that is, for all $n\in \nn$, the set
$$ U^0_{C_n}(M) = \{\p\in\Spec(A_0) \mid \codepth_{(A_0)_\p}(M_\p)\le n\}$$
is open in $\Spec(A_0)$.
\end{proposition*}

\begin{proof} If $A_0$ is a homomorphic image of a regular ring $R_0$,
then the dimension and the depth of the $R_0$-module $M$ are identical to the dimension and depth of $M$
considered as an $R_0$-module. If we show that the set
  $$\widetilde U^0_{C_n}(M)= \{\q\in\Spec(R_0)\mid
  \codepth_{(R_0)_\q}(M_\q)\le n\} $$
  is open in $\Spec(R_0)$ (where
  $M$ is considered a $R_0$-module), then the corresponding set for
  the $A_0$-module $M$ is given by
$$ U^0_{C_n}(M) = \widetilde U^0_{C_n}(M)\cap V(J)$$
where $A_0 = R_0/J$. Thus we may assume that $A_0$ is a regular ring. We may also assume that $A$ is a
polynomial ring over $A_0$ equipped with the standard grading.

Let $\p\in\Spec(A_0)$. By Proposition \ref{AB-formula}, the Auslander-Buchsbaum formula holds:
$$
\depth_{(A_0)_\p}(M_\p) = \depth((A_0)_\p) - \pd_{(A_0)_\p}(M_\p)\,.$$ Let $I= \ann_{A_0}(M)$. By Lemma
\ref{ann}, $I_\p = \ann_{(A_0)_\p}(M_\p)$ and we have that:
$$ \dim_{(A_0)_\p}(M_\p) = \dim((A_0)_\p) - \mbox{ht}(I(A_0)_\p)\,.$$
Suppose that $\p\in\Spec(A_0)$ is such that
$$\codepth_{(A_0)_\p}(M_\p)\le n\,.$$
If $M_\p = 0$ then $\p\not\supseteq I$. Take an element $a\in I\cap (A_0\smallsetminus\p)$. Then for all
$$ \q\in D_a = \{\mathfrak w\in\Spec(A_0)\mid a\notin \mathfrak w\}$$
we have that $M_\q = 0$ and $\codepth_{(A_0)_\q}(M_\q)= -\infty \le n$.

If $M_\p \ne 0$ pick an element $a_1\in A_0\smallsetminus\p$ so that $(A_0)_{\p}$ and $(A_0)_{a_1}$ have the
same minimal primes and put $U_1= D_{a_1} = \{\mathfrak w\in\Spec(A_0)\mid a_1\notin \mathfrak w\}$. Then for
all $\q\in U_1\cap V^0(I)$:
$$\height(I(A_0)_\q) \ge \height(I(A_0)_\p)\,.$$
Let $\pd_{(A_0)_\p}(M_\p)=t$, then by Proposition \ref{Auslander-prop} there is an open subset $U_2$ in
$\Spec(A_0)$ so that
$$\pd_{(A_0)_\q}(M_\q)\le t \quad\text{for all}\quad\q\in U_2\,.$$

Using the Auslander-Buchsbaum formula and the fact that $A_0$ is regular we obtain for all $\q\in U_2\cap
U_1\cap V^0(I)$:
\begin{align*}
\codepth_{(A_0)_\q}(M_\q) &=\dim_{(A_0)_\q}(M_\q) - \depth_{(A_0)_\q}(M_\q) \\ &= \dim((A_0)_\q)
-\mbox{ht}(I(A_0)_\q) - \dim((A_0)_\q) + \pd_{(A_0)_\q}(M_\q) \\ &=\pd_{(A_0)_\q}(M_\q) - \mbox{ht}(I(A_0)_\q).
\end{align*}
This implies that for all $\q\in U = U_1\cap U_2$:
$$ \codepth_{(A_0)_\q}(M_\q) \le \codepth_{(A_0)_\p}(M_\p)$$
and it follows that $U^0_{C_n}(M)$ is an open subset of $\Spec(A_0)$.
\end{proof}

\subsection{A local formula}
\label{local-formula} Using the fact that a complete local Noetherian ring is the homomorphic image of a regular
local ring, we obtain a result similar to \cite[(6.11.5)]{Gr}:

\begin{lemma*} Let $A$ be a Noetherian graded homogeneous ring and let $M$
be a finitely generated graded $A$-module. Then for all prime ideals $\p,\q\in\Spec(A_0)$ with $\p\subseteq \q$
we have that:
$$\codepth_{(A_0)_\q}(M_\q) \ge \codepth_{(A_0)_\p}(M_\p)\,.$$
\end{lemma*}
\noindent {\it Proof.} By replacing $A_0$ by $(A_0)_\q$ (and $A$ by $A_\q$) we may assume that $(A_0, \m_0)$ is
a local ring. Then we have to show:
$$\codepth_{A_0}(M)\ge \codepth_{(A_0)_\p}(M_\p)\,.$$
Let $\widehat\p \in\Spec(\widehat A_0)$ be a minimal prime ideal over $\p\widehat A_0$. In particular,
$\widehat\p\cap A_0 =\p$ and $(\widehat A_0)_{\widehat\p}$ is flat over $(A_0)_\p$ with trivial special fiber.
Moreover:
$$\begin{array}{rcl}
M_\p\otimes_{(A_0)_\p} (\widehat A_0)_{\widehat\p} &= & (\bigoplus_{i\in
\zz} (M_i)_\p) \otimes_{(A_0)_\p} (\widehat A_0)_{\widehat\p} \\
&= &\bigoplus_{i\in \zz} ((M_i)_\p \otimes_{(A_0)_\p} (\widehat
A_0)_{\widehat\p}) \\
&\cong &\bigoplus_{i\in \zz} (\hM_i)_{\widehat\p}
\end{array}
$$
where $\hM_i \cong M_i\otimes_{A_0} \widehat A_0$. We have that:
\begin{gather*}
\depth_{A_0}(M) = \inf\{\depth_{A_0}(M_i)\mid M_i\ne 0\}\\
\dim_{A_0}(M) = \sup\{\dim_{A_0}(M_i)\mid i\in \zz\}
\end{gather*}
By \cite[Theorem 23.3]{Ma}, for all $i\in \zz$:
$$\begin{array}{rcl}
\depth_{(\widehat A_0)_{\widehat\p}}((\hM_i)_{\widehat\p}) &= &\depth_{(A_0)_\p}((M_i)_\p) + \depth((\widehat
A_0)_{\widehat\p}/\p(\widehat A_0)_{\widehat\p}) \\
&=&\depth_{(A_0)_\p}((M_i)_\p)
\end{array}
$$
and by \cite[Theorem 15.1]{Ma}:
$$\begin{array}{rcl}
\dim_{(\widehat A_0)_{\widehat\p}}((\hM_i)_{\widehat\p}) &= &\dim_{(A_0)_\p}((M_i)_\p) + \dim((\widehat
A_0)_{\widehat\p}/\p(\widehat
A_0)_{\widehat\p}) \\
&=&\dim_{(A_0)_\p}((M_i)_\p)
\end{array}
$$
Let
$$ \wM := \bigoplus_{i\in \zz} \hM_i \cong M\otimes_{A_0} \widehat A_0$$
and note that $\wM$ is a finitely generated graded module over the Noetherian homogeneous graded ring
$$\widetilde A := A\otimes_{A_0} \widehat A_0\,.$$
The computation above shows that
$$ \codepth_{(\widehat A_0)_{\widehat\p}} (\wM_{\widehat\p}) =
\codepth_{(A_0)_\p}(M_\p)=:n\,.$$

Since $\widehat A_0$ is a homomorphic image of a regular local ring, by Proposition \ref{dim-formula} the set
$U^0_{C_{n-1}}(\wM)$ is open in $\Spec(\widehat A_0)$. This implies that
$$
\codepth_{\widehat A_0}(\wM) \ge \codepth_{(\widehat
  A_0)_{\widehat\p}}(\wM_{\widehat\p})\,.$$
The same argument as above shows that
$$ \codepth_{\widehat A_0}(\wM)= \codepth_{ A_0}( M)$$
which proves the claim:
\begin{xxalignat}{3}
&\ \ &\codepth_{A_0}(M)\ge&\codepth_{(A_0)_\p}(M_\p)\,. &&\square
\end{xxalignat}

\subsection{Formulas for depth and codepth}
\numberwithin{theorem}{subsection} In this section we make the same assumption as at the beginning, namely, $A$
is a positively graded Noetherian homogeneous ring and $M$ is a finitely generated graded $A$-module. The
following proposition is the graded version of \cite[(6.10.6)]{Gr}:

\begin{proposition}
\label{depth-formula} Let $A$ and $M$ be as above and assume that $A$ is excellent. Then for every
$\p\in\Spec(A_0)$ there is an open subset $U^0\subseteq\Spec(A_0)$ with $\p\in U^0$ so that for all $\q\in
U^0\cap V^0(\p)$:
$$\depth_{(A_0)_\q}(M_\q) = \depth_{(A_0)_\p}(M_\p) +
\depth((A_0)_\q/\p(A_0)_\q)\,.$$
\end{proposition}

\begin{proof} Let $\p\in\Spec(A_0)$, then by Lemma \ref{local-formula} for
all $\q\in V^0(\p)$:
$$ \codepth_{(A_0)_\q}(M_\q) \ge \codepth_{(A_0)_\p}(M_\p),$$
or equivalently:
$$ (^*) \quad \dim_{(A_0)_\q}(M_\q) - \depth_{(A_0)_\q}(M_\q) \ge
\dim_{(A_0)_\p}(M_\p) - \depth_{(A_0)_\p}(M_\p)\,.$$

According to Proposition \ref{dim-formula} there is an open subset $U_1\subseteq\Spec(A_0)$ with $\p\in U_1$ so
that for all $\q\in U_1 \cap V^0(\p)$:
$$ \dim_{(A_0)_\q}(M_\q)= \dim_{(A_0)_\p}(M_\p) + \dim((A_0/\p)_\q)\,.$$

Since $A_0$ is excellent, there is an open subset $U_2\subseteq\Spec(A_0)$ so that $\p\in U_2$ and for all
$\q\in U_2\cap V^0(\p)$ the local ring:
$$(A_0/\p)_\q \quad \mbox{is Cohen-Macaulay.}$$

There is also an open subset $U_3\subseteq\Spec(A_0)$ so that $\p\in U_3$ and for all $\q\in U_3\cap V^0(\p)$ we
have equality on the set of minimal primes:
$$\mbox{Min}_{(A_0)_\q}(I(A_0)_\q) = \mbox{Min}_{(A_0)_\p}(I(A_0)_\p)$$
where $I:=\ann_{A_0}(M)$ denotes the $A_0$-annihilator of $M$. In particular, for all $\q\in U_3\cap V^0(\p)$:
$$ \mbox{ht}(I(A_0)_\q) = \mbox{ht}(I(A_0)_\p)\,.$$

Put $\widetilde U_1 = U_1\cap U_2 \cap U_3$, then for all $\q\in \widetilde U_1\cap V^0(\p)$:
$$ \dim_{(A_0)_\q}(M_\q) = \dim((A_0/I)_\q) \quad \mbox{and} \quad
\dim_{(A_0)_\p}(M_\p) = \dim((A_0/I)_\p).$$

Since $A$ is excellent, the ring $A_0$ is universally catenary and for all $\q\in \widetilde U_1\cap V^0(\p)$:
$$ \dim((A_0/I)_\q) -\dim((A_0/I)_\p) = \dim((A_0/\p)_\q) =
\depth((A_0/\p)_\q)\,.$$

From $(^*)$ we obtain:
$$\depth_{(A_0)_\q}(M_\q) - \depth_{(A_0)_\p}(M_\p)\le
\depth((A_0/\p)_\q)$$ for all $\q\in \widetilde U_1 \cap V^0(\p)$.

In order to prove the other inequality:
$$\depth_{(A_0)_\q}(M_\q) - \depth_{(A_0)_\p}(M_\p)\ge
\depth((A_0/\p)_\q)$$ assume that $\depth_{(A_0)_\p}(M_\p) = t$ and let $f_1,\hdots,f_t\in\p$ be such that
$f_1,\hdots,f_t$ is a regular sequence on $M_\p$. A prime avoidance argument shows that there is an element $a
\in A_0\smallsetminus\p$ so that $f_1,\hdots,f_t$ is a regular sequence on $M_a$. (The argument makes again use
of the fact that the sets $\Ass_{A_0}(M)$ and $\Ass_{A_0}(M/(f_1,\hdots,f_i)M)$ for all $1\le i\le t$ are
finite.)

Put
$$ \ovM := M/(f_1,\dots,f_t)M$$
and consider the associated graded module:
$$ \gr_\p(\ovM) = \bigoplus_{i\in \nn} \p^i\ovM/\p^{i+1}\ovM.$$

The module $\ovM$ is finitely generated over $A$ and $\gr_\p(\ovM)$ is a finitely generated $\gr_\p(A)$-module.
Also note that $\gr_\p(A)$ is a finitely generated algebra over $A/\p A$ and that $A/\p A$ is a finitely
generated algebra over $A_0/\p$. Thus $\gr_\p(A)$ is a finitely generated $A_0/\p$-algebra. By \cite[Theorem
24.1]{Ma} there is an element $b\in A_0 \smallsetminus\p$ so that the $(A_0/\p)_b$-module:
$$ \gr_\p(\ovM)_b = \bigoplus_{i\in \nn} (\p^i\ovM/\p^{i+1}\ovM)_b$$
is free. Set $\widetilde U_2 = D_b = \{\q\in\Spec(A_0) \mid b\notin\q\}$ and fix a prime ideal $\q\in \widetilde
U_2 \cap V^0(\p)$. Assume that
$$\depth((A_0/\p))_\q = s$$
and let $g_1,\hdots, g_s\in \q$ be such that $g_1,\hdots,g_s$ is a regular sequence on $(A_0/\p)_\q$.

\begin{Claim1} $g_1$ is a regular element on $\ovM_\q$.
\end{Claim1}

\begin{Claim2} Set $N_1:= \ovM_\q/g\ovM_\q$, then
$\gr_\p(N_1) \cong \gr_\p(\ovM_\q)/g\gr_\p(\ovM_\q)$.
\end{Claim2}

Assuming the claims, we finish the proof. From the second claim it follows that $\gr_\p(N_1)$ is a free
$(A_0/(g_1,\p)A_0)_\q$-module. Since $g_2$ is a regular element on $(A_0/(g_1,\p)A_0)_\q$, we may apply claims 1
and 2 to $N_1$.  Note that $N_1$ is also a finitely generated graded $A_\q$-module. This yields that $g_2$ is a
regular element on $N_1$ and that with $N_2 = N_1/g_2N_1$:
$$\gr_\p(N_2) \cong \gr_\p(N_1)/g_2\gr_\p(N_1)\,.$$
An induction argument yields that $g_1,\hdots,g_s$ is a regular sequence on $\ovM_\q$ and we have that:
$$\depth_{(A_0)_\q}(M_\q)\ge \depth_{(A_0)_\p}(M_\p) +
\depth((A_0/\p)_\q\,.$$ This inequality holds for all $\q \in \widetilde U_2\cap V^0(\p)$. Assuming the claims
the proposition is now proved with $U^0=\widetilde U_1 \cap \widetilde U_2$.

In order to prove the claims, set $g=g_1$ and $N=N_1$.

\smallskip

\noindent{\it Proof of Claim 1.} Let $z\in \ovM_\q$ with $gz = 0$. Consider the image $\bar z$ of $z$ in
$\ovM_\q/\p\ovM_\q$. Since $\ovM_\q/\p\ovM_\q$ is a free module over $(A_0/\p)_\q$ and since $g$ is regular on
$(A_0/\p)_\q$ we obtain that $\bar z = 0$ and $z\in\p\ovM_\q$. Now consider the image of $z$ in
$\p\ovM_\q/\p^2\ovM_\q$ and repeat the argument. This yields
$$ z\in \bigcap_{j=0}^{\infty} \p^j\ovM_\q\,.$$
Note that
$$\ovM_\q = \bigoplus_{i\in \zz} (\ovM_i)_\q \quad \mbox{with} \quad
(\ovM_i)_\q = (M_i)_\q/(f_1,\hdots,f_t)(M_i)_\q\,.$$ In particular,
$$ \p^j \ovM_\q = \bigoplus_{i\in \zz} \p^j(\ovM_i)_\q$$
and every $ (\ovM_i)_\q$ is a finitely generated $(A_0)_\q$-module. This shows that $z=0$.
\smallskip

\noindent{\it Proof of Claim 2.} By assumption, we have that $\gr_\p(\ovM_\q)$ is a free $(A_0/\p)_\q$- module
and $\p^j\ovM_\q/\p^{j+1}\ovM_\q$ is a direct summand of $\gr_\p(\ovM_\q)$. Thus $\p^j\ovM_\q/\p^{j+1}\ovM_\q$
is a free $(A_0/\p)_\q$-module and $g$ is regular on $(A_0/\p)_\q$. Therefore:
$$ (^{**}) \quad \p^j\ovM_\q \cap g\ovM_\q = g\p^j\ovM_\q$$
and thus:
$$\begin{array}{rcl}
\p^j\ovM_\q/g\p^j\ovM_\q &\cong &\p^j\ovM_\q/ (\p^j\ovM_\q\cap g\ovM_\q)
\\
&\cong & \p^j(\ovM_\q/g\ovM_\q)
\end{array}
$$
From the commutative diagram:
\[
\xymatrixrowsep{0.1pc} \xymatrixcolsep{1.0pc} \xymatrix{
  0\ar@{->}[r]& {\qquad\p^{j+1}N\quad}\ar@{->}[r]\ar@{=}[dd]
&{\quad\p^jN\,\quad\,}\ar@{=}[dd]\ar@{->}[r]&{\qquad\p^jN/\p^{j+1}N\qquad}\ar@{=}[dd]\ar@{->}[r]
&\,0\\
{}&{}&{}&{}&{}&\\
0\ar@{->}[r]& \p^{j+1}(\ovM_\q/g\ovM_\q)\ar@{->}[r]\ar@{->}[ddd]^{\cong}
&\p^j(\ovM_\q/g\ovM_\q)\ar@{->}[ddd]^{\cong}\ar@{->}[r]&\p^j(\ovM_\q/g\ovM_\q)/\p^{j+1}(\ovM_\q/g\ovM_\q)\ar@{->}[ddd]\ar@{->}[r]
&0\\
{}&{}&{}&{}&{}&\\
{}&{}&{}&{}&{}&\\
0\ar@{->}[r]&\p^{j+1}\ovM_\q/g\p^{j+1}\ovM_\q\ar@{->}[r]& \p^j\ovM_\q/g\p^j\ovM_\q
\ar@{->}[r]&\p^j\ovM_\q/(g\p^j\ovM_\q +\p^{j+1}\ovM_\q)\ar@{->}[r] &0 }
\]
we obtain that:
$$\begin{array}{rcl}
\gr_\p(N) & = &\bigoplus_{j\in \nn} \p^jN/\p^{j+1}N \\
&\cong &\bigoplus_{j\in \nn} \p^j\ovM_\q/(g\p^j\ovM_\q + \p^{j+1}\ovM_\q)
\\
&\cong &\bigoplus_{j\in \nn}
(\p^j\ovM_\q/\p^{j+1}\ovM_\q)/g(\p^j\ovM_\q/\p^{j+1}\ovM_\q) \\
&\cong &\text{gr}_\p(\ovM_\q)/g(\gr(\ovM_\q).
\end{array}
$$

This proves the claim, and finishes the proof.
\end{proof}
Similarly to \cite[(6.11.8.1)]{Gr} we have in the graded case:
\begin{corollary}
\label{codepth-formula}
 Let $A$ and $M$ be as above and assume that $A$ is excellent. Then for
every $\p\in\Spec(A_0)$ there is an open subset $U^0\subseteq\Spec(A_0)$ with $\p\in U^0$ so that for all $\q\in
U^0\cap V^0(\p)$:
$$\codepth_{(A_0)_\q}(M_\q) = \codepth_{(A_0)_\p}(M_\p) +
\codepth((A_0)_\q/\p(A_0)_\q)\,.$$
\end{corollary}
\noindent{\it Proof.}Let $\p\in\Spec(A_0)$ and $U_1^0$ be as in Proposition \ref{depth-formula}, so that $\p\in
U_1^0$ and for all $\q\in U_1^0\cap V^0(\p)$:
$$\depth_{(A_0)_\q}(M_\q) = \depth_{(A_0)_\p}(M_\p) +
\depth((A_0)_\q/\p(A_0)_\q)\,.$$

By Proposition \ref{dim-formula} there is an open subset $U^0_2$ in $\Spec(A_0)$ so that $\p\in U^0_2$ and for
all $\q\in U^0_2 \cap V^0(\p)$:
$$ \dim_{(A_0)_\q}(M_\q) = \dim_{(A_0)_\p}(M_\p) + \dim((A_0/\p)_\q)\,.$$

Thus with $U^0 = U^0_1\cap U^0_2$ we have that $\p\in U^0$ and for all $\q\in U^0\cap V^0(\p)$:
\begin{xxalignat}{3}
&\ &\codepth_{(A_0)_\q}(M_\q) &= \codepth_{(A_0)_\p}(M_\p) + \codepth((A_0)_\q/\p(A_0)_\q)\,. &&\square
\end{xxalignat}
We are now ready to prove the graded version of \cite[(6.11.2)(a)]{Gr}:
\begin{theorem}
\label{codepth-locus} Let $A= \bigoplus_{i\in \nn} A_i$ be an excellent graded homogeneous ring and $M =
\bigoplus_{i\in \zz} M_i$ be a finitely generated graded $A$-module. Then for all $n\in \nn$ the set
$$ U^0_{C_n}(M) = \{\p\in\Spec(A_0)\mid \codepth_{(A_0)_\p}(M_\p)\le n
\}$$ is open in $\Spec(A_0)$.
\end{theorem}
\begin{proof} According to Nagata's criterion on openness of loci (see
\cite[Theorem 24.2]{Ma}) we need to show:
\begin{enumerate}[\qquad\,\,\rm(a)]
\item If $\p,\q\in\Spec(A_0)$ with $\q\subseteq \p$ and $\p\in U^0_{C_n}(M)$ then $\q\in U^0_{C_n}(M)$. \item
If $\p\in U^0_{C_n}(M)$ then $U^0_{C_n}(M)$ contains a nonempty open subset of $V(\p)$.
\end{enumerate}

$($a$)$ Let $\p,\q\in\Spec(A_0)$ with $\q\subseteq \p$. By Lemma \ref{local-formula}:
$$\codepth_{(A_0)_\p}(M_\p) \ge \codepth_{(A_0)_\q}(M_\q)$$
and thus $\p\in U^0_{C_n}(M)$ implies that $\q\in U^0_{C_n}(M)$.

$($b$)$ Let $\p\in U^0_{C_n}(M)$. By Corollary \ref{codepth-formula} there is an open subset $U^0_1$ in
$\Spec(A_0)$ so that $\p\in U^0_1$ and for all $\q\in U^0_1\cap V^0(\p)$:
$$\codepth_{(A_0)_\q}(M_\q) = \codepth_{(A_0)_\p}(M_\p) +
\codepth((A_0)_\q/\p(A_0)_\q)\,.$$ Since $A$ and $A_0$ are excellent, there is an open subset $U^0_2$ in
$\Spec(A_0)$ so that $\p\in U^0_2$ and for all $\q\in U^0_2\cap V^0(\p)$ the ring $(A_0/\p)_\q$ is
Cohen-Macaulay. Therefore with $U^0 = U^0_1\cap U^0_2$ we have that $\p\in U^0$ and for all $\q\in U^0\cap
V^0(\p)$:
$$\codepth_{(A_0)_\q}(M_\q) = \codepth_{(A_0)_\p}(M_\p)\,.$$
This implies that $U^0\cap V^0(\p)\subseteq U^0_{C_n}(M)$ and the Theorem is proved.
\end{proof}
\begin{corollary}
\label{CM-locus} Let $A$ and $M$ be as in Theorem {\rm\ref{codepth-locus}}. Then the Cohen-Macaulay locus of the
$A_0$-module $M$:
$$ U^0_{CM}(M) =U^0_{C_0}(M) = \{\p\in\Spec(A_0)\mid M_\p \; \mbox{is a CM module
over}\; (A_0)_\p\}$$ is open in $\Spec(A_0)$.\qed
\end{corollary}

\section{Openness of the $(S_n)$-locus}
\numberwithin{theorem}{section} Throughout this section we assume that $R=A_0$ is the base ring of a graded
Noetherian homogeneous ring $A= \bigoplus_{i\ge 0} A_i$ and $M$ is a finitely generated graded $A$-module. This
includes the case of a finitely generated module $M$ over a Noetherian ring $R$. For those modules we prove that
the openness of the $C_n$-loci of $M$ implies the openness of the $(S_k)$-loci of $M$. The argument is due to
Grothendieck \cite[(5.7.2) and (6.11.2)(b)]{Gr} but we include it here for the convenience of the reader. The
proof also shows that the $(S_k)$-loci of $M$ only depend on the $C_n$-loci of $M$ and on the annihilator of
$M$, so that two $R$-modules $M$ and $N$ with the same annihilators and $C_n$-loci have identical $(S_k)$-loci.

Let $M$ be an $R$-module and suppose that for all $n\in \nn_0$ the set:
$$U_{C_n}(M) = \{\p\in\Spec(R)\mid \codepth_{R_\p}(M_\p)\le n \}$$
is open in $\Spec(R)$. Define:
$$ Z_n = V(\fb_n) =\Spec(R)\smallsetminus U_{C_n}(M)$$
where $\fb_n\subseteq R$ is a reduced ideal. Obviously, for all $n\in \nn$:
$$U_{C_n}(M)\subseteq U_{C_{n+1}}(M)$$
and therefore:
$$Z_{n+1}\subseteq Z_n \quad \mbox{and} \quad \fb_n \subseteq
\fb_{n+1}\,.$$ Since $R$ is Noetherian there is an $m\in \nn$ so that for all $t\in \nn$:
$$\fb_m = \fb_{m+t} \quad \mbox{and} \quad Z_m = Z_{m+t}\,.$$

\begin{lemma}
\label{Z}Let $m\in \nn$ be as above. Then $Z_m =\emptyset$.
\end{lemma}
\begin{proof} If $\p\in Z_m$ then $\p\in Z_{m+t}$ for all $t\in \nn$. By
definition of $Z_{n+t}$:
$$ \codepth_{(R)_\p}(M_\p)\ge m+t \quad \mbox{for all} \quad t\in
\nn\,.$$ But $\codepth_{R_\p}(M_\p)\le \dim((R)_\p)\le \infty$ and therefore $Z_m=\emptyset$.
\end{proof}
Recall that the $R$-module $M$ satisfies Serre's condition $(S_k)$ if for all $\p\in \Spec(R)$:
$$ (^*) \quad \depth_{R_\p}(M_\p) \ge \mbox{min}(\mbox{dim}(M_\p),k).$$
From now on let $m$ denote the minimal $m\in \nn$ with $Z_m =\emptyset$.

\begin{lemma}
\label{b} With the assumptions as above put $\ov R = R/\ann_R(M)$ and let $k\in \nn$. Then the $R$-module $M$
satisfies $(S_k)$ if and only if for all $0\le n<m$:
$$\height(\fb_n\ov R) > n+k\,.$$
\end{lemma}
\begin{proof}
Suppose that $M$ satisfies $(S_k)$ and fix an integer $n$ with $0\le n<m$. Let $\p\in\Spec(R)$ with
$\fb_n\subseteq \p$. Then $\p\in Z_n$ and therefore:
$$\codepth_{R_\p}(M_\p)>n$$
or equivalently:
$$\dim_{R_\p}(M_\p)- \depth_{R_\p}(M_\p) >n\,.$$
Since $M$ satisfies $(S_k)$ we obtain that whenever
$$\dim_{R_\p}(M_\p)- \depth_{R_\p}(M_\p)\ne 0$$
then
$$\depth_{R_\p}(M_\p)\ge k\,.$$
Thus, if $\p\in Z_n$ then
$$\dim_{R_\p}(M_\p)\ge n+k$$
which implies that $\mbox{ht}(\fb_n\ov R)\ge n+k$.

Conversely, fix an integer $k$ and assume that for all $0\le n<m$:
$$\mbox{ht}(\fb_n\ov R)>n+k\,.$$
Let $\p\in\Spec(R)$.

If $M_\p = 0$ then $\depth_{R_\p}(M_\p) = \infty$ and condition $(^*)$ is satisfied.

Now assume $M_\p\ne 0$. If $M_\p$ is a Cohen-Macaulay $R$-module, then condition $(^*)$ is satisfied. Now assume
that:
$$\codepth_{R_\p}(M_\p)>0$$ and let $n\in\nn_0$ with
$$\codepth_{R_\p}(M_\p) = n+1\,.$$
Thus $\p\in Z_n$ and $\fb_n\subseteq \p$. By assumption
$$\mbox{ht}(\fb_n\ov R)>n+k \; \Rightarrow \;
\mbox{ht}(\fb_nR_\p)>n+k \; \Rightarrow \; \mbox{dim}(\ov R_\p) > n+k\,.$$ This implies that
$$\begin{array}{rcl}
\codepth_{R_\p}(M_\p) &= &n+1 \\
&= &\dim(\ov R_\p) - \depth_{R_\p}(M_\p) \\
&\ge &n+1+k - \depth_{R_\p}(M_\p)
\end{array}
$$
and therefore
$$ \depth_{R_\p}(M_\p)\ge k\,.$$
Thus $M_\p$ satisfies condition $(^*)$ and the $R$-module $M$ satisfies Serre's condition $(S_k)$.
\end{proof}

For all $0\le n<m$ consider the closed subset of $\Spec(R)$:
$$Y_{n,k} = \{ \q\in V(\fb_n) \mid \mbox{ht}(\fb_n\ov R_\q) \le n+k\}\,, $$
and its complement
 $$V_{n,k} = \Spec(R) -Y_{n,k}\,,$$ an open subset of $\Spec(R)$. By Lemma \ref{b}:
$$U_{S_k}(M) = \bigcap_{0\le n<m} V_{n,k}\,$$
is an open subset of $\Spec(R)$. We have shown:
\begin{theorem}
\label{S-locus} Let $M$ be an $R$-module as above. If for all $n\in \nn_0$ the $C_n$-locus $U_{C_n}(M)$ is open
in $\Spec(R)$ then for all $k\in \nn$ the $(S_k)$-locus:
$$U_{S_k}(M) = \{ \p\in \Spec(R) \mid M_\p \; \mbox{satisfies} \; (S_k)\}$$
is open in $\Spec(R)$.\end{theorem} In the graded case the theorem states:
\begin{corollary}
\label{S-locus2} Let $A= \bigoplus_{i\in \nn} A_i$ be an excellent graded homogeneous ring and let $M =
\bigoplus_{i\in \zz} M_i$ be a finitely generated graded $A$-module. Then for all $k\in \nn$ the set
$$ U^0_{S_k}(M) = \{\p\in\Spec(A_0)\mid \; \mbox{the $(A_0)_\p$-module}
\; M_\p \; \mbox{satisfies} \; (S_k) \}$$ is open in $\Spec(A_0)$.
\end{corollary}
The proof of the theorem also yields the following corollary:
\begin{corollary}\label{same-modules} Suppose that $M$ and $N$ are $R$-modules as above. Assume that $\mbox{ann}_R(M) =
\mbox{ann}_R(N)$ and that for all $n\in \nn_0$ the sets $U_{C_n}(M)=U_{C_n}(N)$ are open in $\Spec(R)$. Then for
all $k\in \nn$:
$$ U_{S_k}(M) = U_{S_k}(N)$$
and the $(S_k)$-loci are open subsets of $\Spec(R)$. \end{corollary}
\section{Stability on the homogeneous parts}
\numberwithin{theorem}{section} Let $ A = \bigoplus_{i\in \nn} A_i$ be an excellent graded homogeneous
Noetherian ring and let $M = \bigoplus_{i\in \zz} M_i$ be a finitely generated graded $A$-module. In this
section we prove that there is a $k\in \nn$ so that for all $n\in \nn$ and all $i\ge k$:
$$ U^0_{C_n}(M_i) = U^0_{C_n}(M_k) \quad \mbox{and} \quad U^0_{S_n}(M_i) =
U^0_{S_n}(M_k),$$ that is, the codepth and $(S_n)$ loci of the homogeneous parts of $M$ are eventually stable
(considered as an $A_0$-module). As before we define for all $t\in \zz$:
$$ N_t = \bigoplus_{i\ge t} M_i,$$
and observe the following simple facts: Let $k_1\in \nn$ be an integer so that for all $t\ge k_1$:
$\ann_{A_0}(M_t) = \ann_{A_0}(M_{k_1})$. Then for all $t\ge k_1$:
$$U^0_{C_n}(N_t) \supseteq U^0_{C_n}(N_{k_1}) \quad \mbox{and} \quad
U^0_{S_n}(N_t) \supseteq U^0_{S_n}(N_{k_1})\,.$$ Since $A_0$ is Noetherian there is an integer $k_2\in \zz$ so
that $k_2\ge k_1$ and
$$U^0_{C_n}(N_t)= U^0_{C_n}(N_{k_2}) \quad \mbox{and} \quad
U^0_{S_n}(N_t)= U^0_{S_n}(N_{k_2})\,.$$ We may also assume for large enough $k_2$ that:
$$N_{k_2} = A M_{k_2},$$
which implies that for all $t\ge k_2$:
$$N_t = A M_t.$$

\begin{lemma}
\label{stable-depth} With the assumptions as above assume additionally that $(A_0,\m_0)$ is a local ring. Then
there is a $k_3\in \zz$ so that for all $t\ge k_3$:
$$\depth_{A_0}(M_t) = \depth_{A_0}(M_{k_3}) =\depth_{A_0}(N_{k_3})\,.$$
\end{lemma}

\noindent {\it Proof.} Let $k_1$ and $k_2$ be as above and take an integer $k$ with $k>k_2$. Then
$\codepth_{A_0}(N_k) =n$ for some $n\in \nn$ and therefore:
$$\m_0\in U^0_{C_n}(N_k) \quad \mbox{and} \quad \m_0\notin
U^0_{C_{n-1}}(N_k)\,.$$ Since $k\ge k_2$ we have for all $t\ge k$:
$$ \codepth_{A_0}(N_k) = n = \codepth_{A_0}(N_t)\,.$$
For all $t\ge k_1$ we also have that $\ann_{A_0}(N_t) = \ann_{A_0}(N_k)$, and therefore for all $t\ge k$:
$$ \depth_{A_0}(N_t) =s = \depth_{A_0}(N_k)\,.$$

Let $r_1,\hdots,r_s$ be a maximal regular sequence on $N_k$ and put
$$ \ovN_k= N_k/(r_1,\hdots,r_s)N_k \quad \mbox{with homogeneous parts}
\quad \ovM_i = M_i/(r_1,\hdots,r_s)M_i$$ for $i\ge k$. Note that the torsion submodule $\Gamma_{A_+}(\ovN_k)$ is
a finitely generated $A$-submodule of $\ovN_k$. This implies that there is an integer $k_3\ge k$ so that
$\Gamma_{A_+}(\ovN_k)\cap N_{k_3} = 0 = \Gamma_{A_+}(\ovN_{k_3})$. Thus for $k_3$ large enough the $A$-module
$\ovN_{k_3}$ is $A_+$-torsion free. Since by assumption $\depth_{A_0}(N_k) =s = \depth_{A_0}(N_{k_3})$ there is
an integer $i\ge k_3$ and an element $\bar x\in \ovM_i$ so that $\bar x\ne 0$ and $\m_0\bar x = 0$. Since
$\ovN_{k_3}$ is $A_+$-torsion free we obtain:
$$(A_+)^l\bar x \ne 0 \quad \mbox{for all} \quad l\in \nn\,.$$
Thus for $k_4=i>k_3$ we have that $\depth_{A_0}(\ovM_{k_4+l}) = 0$ for all $l\in \nn_0$ and therefore for all
$t\ge k_4$:
\begin{xxalignat}{3}
&\ &\depth_{A_0}( M_t)&=\depth_{A_0}(M_{k_4}) =s  &&\square
\end{xxalignat}

Choose an integer $k_0\in \zz$ so that the following conditions are satisfied:
\begin{enumerate}[\qquad \rm (a)]
\item $N_{k_0} = AM_{k_0}$, that is, $N_{k_0}$ is generated in the lowest nonvanishing degree. \item For all
$t\ge k_0$: $\ann(M_{k_0}) = \ann(M_t)$. \item For all $n\in \nn_0$ and all $t\ge k_0$:
$$U^0_{C_n}(N_t) = U^0_{C_n}(N_{k_0}) \quad \mbox{and} \quad U^0_{S_n}(N_t) =
U^0_{S_n}(N_{k_0})\,.$$
\end{enumerate}
As before put:
$$ Z_n =\Spec(A_0)\smallsetminus U^0_{C_n}(N_{k_0})  = V(\fb_n)$$
where $\fb_n\subseteq A_0$ is a reduced ideal. Then $\fb_n \subseteq \fb_{n+1}$ yielding an increasing sequence
of ideals:
$$ \fb_0\subseteq \fb_1 \subseteq \hdots \subseteq \fb_{m-1} \subseteq
\hdots$$ We have seen before that the sequence stops with some $\fb_m = A_0$ and let $m$ be minimal with this
property, that is, let $\fb_m= A_0$ and $\fb_{m-1}\ne A_0$. For all $0\le j\le m-1$ we consider the set of
minimal prime divisors of $\fb_j$:
$$ \mbox{Min}(A_0/\fb_j) = \{\p_{j1},\hdots,\p_{jr_j}\}\,.$$
By Lemma \ref{stable-depth} for all $0\le j\le m-1$ and all $r_j\ge h\ge 1$ there is an integer $k_{jh}\in \nn$
with $k_{jh}\ge k_0$ so that for all $i\ge k_{jh}$:
$$ \depth_{(A_0)_{\p_{jh}}}((M_i)_{\p_{jh}} =
\depth_{(A_0)_{\p_{jh}}}((M_{k_{jh}})_{\p_{jh}} =\mbox{constant}\,.$$ Let $k=\mbox{max}\{k_{jh}\mid 0\le j\le
m-1; 1\le h\le r_j\}$, then for all $i \ge k$:
$$ \depth_{(A_0)_{\p_{jh}}}((M_i)_{\p_{jh}}) =
\depth_{(A_0)_{\p_{jh}}}((M_k)_{\p_{jh}})=\depth_{(A_0)_{\p_{jh}}}((N_k)_{\p_{jh}})\,.$$ By assumption on the
annihilators we also have for all $i \ge k$:
$$\dim_{(A_0)_{\p_{jh}}}((M_i)_{\p_{jh}}) =
\dim_{(A_0)_{\p_{jh}}}((M_k)_{\p_{jh}})=\dim_{(A_0)_{\p_{jh}}}((N_k)_{\p_{jh}})$$ which implies that for all
$i\ge k$ and all primes $\p_{jh}$:
$$\codepth_{(A_0)_{\p_{jh}}}((M_i)_{\p_{jh}}) =
\codepth_{(A_0)_{\p_{jh}}}((M_k)_{\p_{jh}})=\codepth_{(A_0)_{\p_{jh}}}((N_k)_{\p_{jh}})\,.$$ We are now ready to
prove:

\begin{theorem}
\label{stable-codepth}
 Let $k$ be as above. Then for all $i\ge k$ and all $\p\in\Spec(A_0)$:
$$\codepth_{(A_0)_\p}((M_i)_\p) = \codepth_{(A_0)_\p}((M_k)_\p)\,.$$
\end{theorem}

\noindent {\it Proof.} Let $\p\in\Spec(A_0)$. If $\fb_0\not\subseteq \p$ then $(N_k)_\p$ is a Cohen-Macaulay
module over$(A_0)_\p$ it follows that $(M_i)_\p$ is Cohen-Macaulay for all $i\ge k$.

Assume that $\fb_0\subseteq \p$ and let $g$ be minimal so that $\fb_g\subseteq \p$ and $\fb_{g+1}\not\subseteq
\p$. In this case $\codepth_{(A_0)_\p}((N_k)_\p) = g+1$ and there is an integer $1\le j\le r_j$ so that
$\p_{gj}\subseteq \p$. By  \cite[(6.11.5)]{Gr}, the nongraded version of Lemma \ref{local-formula}, for all
$i\ge k$:
$$ \codepth_ {(A_0)_\p}((M_i)_\p) \ge \codepth_
{(A_0)_{p_{gj}}}((M_i)_{p_{gj}})=\depth_{(A_0)_{\p_{gj}}}((N_k)_{\p_{gj}})> g\,.$$ In order to verify the other
inequality consider
$$ \codepth_{(A_0)_\p}((N_k)_\p) = g+1 = \dim((N_k)_\p) -
\depth_{(A_0)_\p}((N_k)_\p)$$ and assume that $\depth_{(A_0)_\p}((N_k)_\p) =s$. Let $x_1,\hdots,x_s$ be a
regular sequence on $(N_k)_\p$, then  $x_1,\hdots,x_s$ is a regular sequence on $(M_i)_\p$ for all $i\ge k$.
Since $N_k$ and $M_i$ have the same annihilators we obtain that:
$$\codepth_{(A_0)_\p}((N_k)_\p) =g+1 \ge \codepth_{(A_0)_\p}((M_i)_\p)$$
for all $i\ge k$.  This shows that for all $i\ge k$:
\begin{xxalignat}{3}
&\ &\codepth_{(A_0)_\p}((M_i)_\p)&= g+1\,.&&\square
\end{xxalignat}

\begin{corollary} There is an integer $k\in \nn$ so that for all $i\ge k$
and all $n\in \nn$:
\begin{xxalignat}{3}
&\ & U^0_{C_n}(M_i) = U^0_{C_n}(M_k)& = U^0_{C_n}(N_k) \,. && \square
\end{xxalignat}
\end{corollary}

\begin{corollary} There is an integer $k\in \nn$ so that for all $i\ge k$
and all $n\in \nn$:
$$ U^0_{S_n}(M_i) = U^0_{S_n}(M_k) = U^0_{S_n}(N_k)\,.$$
\end{corollary}

\begin{proof} The second corollary follows from the first by using Corollary {\rm\ref{same-modules}}.
\end{proof}

\section{Applications}

Let $A$ be an excellent ring, $M$ a finitely generated $A$-module, and let $I\subseteq A$ be an ideal of $A$. By
applying the results of the previous section to the Rees algebra/module and to the associated graded
ring/module, respectively, we see that there is an integer $k\in \nn$ so that for all $i\ge k$ and all $n\in
\nn$:
\begin{align*}
&U_{C_n}(I^iM) = U_{C_n}(I^kM) \qquad\text{and}  &U_{C_n}(I^iM/I^{i+1}M) =
U_{C_n}(I^kM/I^{k+1}M)\\
&U_{S_n}(I^iM) = U_{S_n}(I^kM) \,\qquad\text {and} &U_{S_n}(I^iM/I^{i+1}M) = U_{S_n}(I^kM/I^{k+1}M)
\end{align*}

In the following we want to apply these results to the $(S_n)$- and codepth-loci of the modules $M/I^kM$. We
want to show that these loci are again eventually stable provided that $M$ is a Cohen-Macaulay module over $A$.
\begin{lemma}\label{supp} Let $A$ be any Noetherian ring, $I\subseteq A$ an ideal, and
$M$ a finitely generated $A$-module. Then for all $k\in \nn$:
$$ \Supp(M/I^kM) = \Supp(M/IM)\,.$$
\end{lemma}
\begin{proof} It suffices to show that for all $k\in \nn$:
$$ \Supp(M/I^kM) = \Supp(M/I^{k+1}M)\,.$$
Since $M/I^kM$ is a homomorphic image of $M/I^{k+1}M$ we have $\Supp(M/I^kM)\subseteq  \Supp(M/I^{k+1}M)$.
Consider the exact sequence:
$$ 0\to I^kM/I^{k+1}M \to M/I^{k+1}M \to M/I^kM \to 0$$ and let
$\p\in\Spec(A)$ with $I\subseteq \p$. The sequence stays exact after localization:
$$ 0\to (I^kM/I^{k+1}M)_\p \to (M/I^{k+1}M)_\p \to (M/I^kM)_\p \to 0\,.$$
If $(M/I^kM)_\p = 0$ with $(M/I^{k+1}M)_\p\ne 0$, then
$$ (I^kM/I^{k+1}M)_\p = (M/I^{k+1}M)_\p,$$
which implies by Nakayama that $(M/I^{k+1}M)_\p =0$, a contradiction.
\end{proof}

A more general version of the next result was proved, using different methods, by Kodiyalam \cite[Corollary
9]{Ko}.

\begin{theorem} Suppose that $(A,\m)$ is a local Noetherian ring,
$I\subseteq A$ an ideal of $A$, and let $M$ be a finitely generated $A$-module. Then there is a $k\in \nn$ so
that for all $i\ge k$:
  $$\depth_A(M/I^iM) = \depth_A(M/I^kM)\,.$$
\end{theorem}
\noindent{\it Proof.} Let $\widehat A$ be the $\m$-adic completion of $A$. Then for any finitely generated
$A$-module $T$:
$$\depth_A(T) = \depth_{\widehat A}(T\otimes_A \widehat A)$$
and we may replace $A$ by $\widehat A$ and $M$ by $M\otimes_A \widehat A$ and assume that $A$ is excellent. By
Lemma \ref{stable-depth} there is a $k_1\in \nn$ so that for all $t\ge k_1$:
$$\depth_A(I^tM/I^{t+1}M) = \depth_A(I^{k_1}M/I^{k_1+1}M)=g.$$
For all $t\ge k_1$ consider the exact sequence:
$$0 \to I^tM/I^{t+1}M \to M/I^{t+1}M \to M/I^tM \to 0$$
which leads to an exact sequence on the cohomology modules:
\begin{gather*}\cdots\to H^i_\m(M/I^{t+1}M) \to H^i_\m(M/I^tM) \to 0 \to
\cdots \to 0 \to \cdots \\\cdots \to H^{g-1}_\m(M/I^{t+1}M)\to H^{g-1}_\m(M/I^tM)\to H^g_\m(I^tM/I^{t+1}M)\to\\
\to H^g_\m(M/I^{t+1}M) \to H^g_\m(M/I^tM) \to \cdots \end{gather*} where $g$ is minimal with
$H^g_\m(I^tM/I^{t+1}M)\ne 0$.

{\it case 1:} There is an $i\leq g-1$ and a $t_0\ge k_1$ so that $H^i_\m(M/I^{t_0}M)\ne 0$. Then for all $t\ge
t_0$ $H^i_\m(M/I^tM)\ne 0$. Let $h\le g-1$ be the minimal $i$ with this property, then
$$\depth_A(M/I^tM) = h \quad \mbox{for all} \quad t\ge t_0\,.$$

{\it case 2:}  For all $i\le g-1$ and all $t\ge k_1$:
$$H^i_\m(M/I^tM)=0\,.$$
This implies that $\depth_A(M/I^tM) \ge g-1$ for all $t\ge k_1$.

{\it case 2.1:}  There are infinitely many $t\ge k_1$ so that
$$H^{g-1}_\m(M/I^tM)\ne 0\,.$$
From the long exact sequence we observe that $H^{g-1}_\m(M/I^tM)\ne 0$ implies that $H^{g-1}_\m(M/I^{t-1}M)\ne
0$ whenever $t-1\ge k_1$. Thus in this case there is a $t_1\ge k_1$ so that for all $t\ge t_1$:
$$H^{g-1}_\m(M/I^tM)\ne 0\,,$$
and therefore for all $t\ge t_1$: $\depth_A(M/I^tM) = g-1$.

{\it case 2.2:} There is a $t_2\ge k_1$ so that for all $t\ge t_2$: $H^{g-1}_\m(M/I^tM)=0$. Then for all $t\ge
t_2$:
\begin{xxalignat}{3}
&\ &\depth_A(M/I^tM)& = g\,. &&\square
\end{xxalignat}

\begin{theorem} Let $A$ be an excellent ring and $M$ a finitely generated
Cohen-Macaulay $A$-module. Let $I\subseteq A$ be an ideal of $A$ which is not contained in any minimal prime
ideal of $M$. Then there is an integer $k\in \nn$ so that for all $t\ge k$ and all $n\in \nn_0$:
\begin{enumerate}[\quad\rm(1)]
\item $U_{C_n}(M/I^tM) = U_{C_n}(M/I^{k_0}M)$. \item $U_{S_n}(M/I^tM) = U_{S_n}(M/I^{k_0}M)$.
\end{enumerate}
\end{theorem}
\begin{proof} $($1$)$ Fix $n\in \nn$ and let $k\in \nn$ so that for all
$t\ge k$:
$$U_{C_n}(I^tM) = U_{C_n}(I^kM)\,.$$
We claim that for all $i\ge k$ and all $\p\in V(I)$:
$$\depth_{A_\p}((M/I^iM)_\p) = \depth_{A_\p}((M/I^kM)_\p)\,.$$

Obviously, for all $i\ge k$: $\dim((I^iM)_\p) = \dim((I^kM)_\p)$ and thus because of the stability of the
codepth-loci we have for all $\p\in V(I)$ and all $i\ge k$ that:
$$\depth_{A_\p}((I^iM)_\p) = \depth_{A_\p}((I^kM)_\p)\,.$$
Fix an integer $i\ge k$ and a prime ideal $\p\in V(I)$ and consider the exact sequence:
$$0\to (I^iM)_\p \to M_\p \to (M/I^iM)_\p  \to 0\,.$$
With $d=\dim_{A_\p}(M_\p) = \depth_{A_\p}(M_\p)$ we obtain a long exact sequence of the local cohomology
modules:
\begin{gather*} \cdots \to 0\to H^{i-1}_\p((M/I^iM)_\p)\to
H^i_\p((I^iM)_\p) \to 0 \to \cdots \to 0\to \\
\to H^{d-1}_\p((M/I^iM)_\p)\to H^d_\p((I^iM)_\p)\to H^d_\p(M_\p)\to  0 = H^d_\p((M/I^iM)_\p)
\end{gather*}
where $H^d_\p((M/I^iM)_\p)=0$ since $\dim_{A_\p}((M/I^iM)_\p)\le d-1$. This shows that
$$\depth_{A_\p}((M/I^iM)_\p) = \depth_{A_\p}((I^iM)_\p)-1
=\depth_{A_\p}((I^kM)_\p)-1$$ and the claim is proven. For all $i\ge k$ and all $\p\in V(I)$ we have:
\begin{gather*}
\depth_{A_\p}((M/I^iM)_\p)=\depth_{A_\p}((M/I^kM)_\p)\\
\dim((M/I^iM)_\p)= \dim((M/I^kM)_\p)\,.
\end{gather*}
The last equation is obtained from Lemma \ref{supp}. This yields that for all $n\in \nn$ and for all $i\ge k$:
$$U_{C_n}(M/I^iM) = U_{C_n}(M/I^kM)\,.$$

The second assumption follows with Corollary {\rm\ref{same-modules}}.
\end{proof}

\begin{corollary} Let $A$, $M$, and $I$ be as in the theorem and assume that $IM\ne M$. Then there
is an element $a\in A$ so that for all $k\in \nn$:
\begin{enumerate}[\quad\rm(1)]
\item  $(M/I^kM)_a \ne 0$. \item $(M/I^kM)_a$ is a Cohen-Macaulay module.\qed
\end{enumerate}
\end{corollary}
\begin{corollary} Let $A$ be an excellent ring and $M$ a finitely generated $A$-module. Suppose that the ideal
$I\subseteq A$ satisfies the following conditions: \begin{enumerate}[(i)] \item $I$ is not contained in a
minimal prime of $M$. \item If $\frak a \subseteq A$ is the defining ideal of the non-Cohen-Macaulay locus of
$M$ then $\frak a \nsubseteq \sqrt{(IM:M)}$.\end{enumerate} Then there is an element $a\in A$ so that for all
$k\in \nn$:
\begin{enumerate}[\quad\rm(1)]
\item  $(M/I^kM)_a \ne 0$. \item $(M/I^kM)_a$ is a Cohen-Macaulay module.\end{enumerate} \end{corollary}
\begin{proof} Choose an element $b\in \frak a \smallsetminus \sqrt{(IM:M)}$. In order to prove the assertion apply the
previous corollary to the Cohen-Macaulay $A_b$-module $M_b$. \end{proof}

\end{document}